\documentclass[a4paper]{article}
\usepackage{amsmath, amsthm, amsfonts,bm}
\usepackage{fullpage}
\usepackage{bm}
\usepackage{xcolor}
\usepackage{amssymb}
\usepackage{thm-restate}
\usepackage{graphicx}
\usepackage{lineno}
\usepackage{enumitem}
\usepackage{longtable}
\usepackage{environ}
\usepackage{dashrule} 
\usepackage{bm}
\usepackage{thm-restate}
\usepackage{relsize}
\usepackage{enumitem}
\usepackage{tikz}
\usetikzlibrary{shapes.geometric}
\usepackage[numbers,sort&compress]{natbib}
\usepackage{hyperref}
\usepackage{cleveref}
\usepackage{xurl}

\hypersetup{
	colorlinks=true, 
	linkcolor=cyan,
	citecolor=blue,
}

\setlist{  
	listparindent=1.5em}

\newtheorem{thm}{Theorem}
\newtheorem{theorem}{Theorem}[section]

\newtheorem{corollary}[thm]{Corollary}
\newtheorem{problem}[theorem]{Problem} 
\newtheorem{question}[problem]{Question}

\newtheorem{conjecture}[thm]{Conjecture}

\theoremstyle{definition}

\theoremstyle{remark}

\newcommand{\dist}{\operatorname{dist}}

\newcommand{\im}{\operatorname{Im}}
\newcommand{\C}{$\mathcal{C}_4$-free }
\newcommand{\Cc}{$\mathcal{C}_4$-free}
\newcommand{\Ccc}{$\mathcal{C}^*_4$-free}

\newcommand{\s}{\diamond}

  \newcommand{\td}[1]{\hat{#1}}

\title{Oriented Trees in Digraphs without Oriented $4$-cycles
}
\author{Maya Stein\thanks{Department of Mathematical Engineering and Center for Mathematical Modeling (CNRS IRL2807), University of Chile. Supported by FONDECYT Regular Grant 1221905,  by ANID Basal Grant CMM FB210005,  and by MSCA-RISE-2020-101007705 project {\it RandNET}. {\tt mstein@dim.uchile.cl}
	} \and
	Ana Trujillo-Negrete\thanks{Center for Mathematical Modeling (CNRS IRL 2807), University of Chile, Santiago, Chile.  Supported by ANID/Fondecyt Postdoctorado 3220838, and by ANID Basal Grant CMM FB210005. {\tt ltrujillo@dim.uchile.cl } }} 

\date{}
\begin{document}
	\maketitle

	\abstract{We prove that if $D$ is a digraph  of maximum outdegree and indegree at least $k$, 
		and minimum semidegree at least $k/2$ that contains no oriented $4$-cycles, then $D$  contains each oriented tree $T$ with~$k$ arcs. 
		This can be slightly improved if $T$ is either antidirected or an arborescence.
	}

%%%%%%%%%%%%%%%%%%%%%%
%%%%%%%%%%%%%%%%%%%%%
%%%%%%%%%%%%%%%%%%%%
%%%%%%%%%%%%%%%%%%%
%%%%%%%%%%%%%%%%%%
%%%%%%%%%%%%%%%%%
%%%%%%%%%%%%%%%%
%%%%%%%%%%%%%%%
%%%%%%%%%%%%%%
%%%%%%%%%%%%%
%%%%%%%%%%%%
%%%%%%%%%%%
%%%%%%%%%%
%%%%%%%%%
%%%%%%%%
%%%%%%%
%%%%%%
%%%%%
%%%%
%%%
%%
%

%
%
\section{Introduction} 
 Any undirected graph $G$ of minimum degree $\delta(G)$ at least $k$ contains each $k$-edge tree: Indeed, we can just greedily embed the tree. Because of the $k$-edge stars, the bound on $\delta(G)$ cannot be made lower than $k$, unless  some extra condition on $G$ is added. Clearly, one necessary condition is that $G$ has a vertex that can accommodate the highest degree vertex of $T$, in other words, $G$ has maximum degree $\Delta (G)\ge \Delta(T)$. This alone is not enough, which can be seen for instance by considering  two copies of the clique  on $\lfloor k/2\rfloor-1$ vertices and adding a universal vertex. This graph does not contain any $k$-edge tree having a vertex with three   roughly equal-sized components attached to it. Therefore, we should either change  the bounds on $\Delta(G)$ and/or $\delta(G)$ (for a discussion of this option see Section~\ref{sec:otherdegree}), or add more conditions on $G$ apart from  $\delta (G)\ge k/2$ and $\Delta (G)\ge \Delta(T)$. 

A natural candidate for such an additional condition is to ask for some kind of expansion of~$G$. In particular, assuming that $G$ has large girth should be helpful for finding trees. 
 Brandt and Dobson~\cite{erdos-sos-girth5} showed that every graph of girth at least five, with $\delta(G)\ge k/2$ and with  $\Delta(G)\ge \Delta (T)$ contains each $k$-edge tree. 
Generalising 
this result, 
 Sacl\'e and Wo\'zniak~\cite{sacle} proved\footnote{The main result of~\cite{sacle} is that every graph $G$ that has more than $\frac{k-1}2|V(G)|$ edges and no $C_4$-cycles  contains each $k$-edge tree. The result we are referring to here is alluded to in Section~4 of~\cite{sacle}.} 
 the following: If $T$ is a tree with $k$ edges, and $G$ is a graph with $\delta(G)\ge k/2$ and $\Delta(G)> \Delta(T)$ (or $\Delta(G)\ge k$ if $T$ is a star)
  such that~$G$ contains no $4$-cycles, then $T\subseteq G$. 
Extending ideas from~\cite{sacle}, we prove an analogous result for digraphs. 

The \textit{minimum semidegree} $\delta^0(D)$ of a digraph $D$
  is the minimum of the indegrees and outdegrees of the vertices of $D$, and its {\it maximum total degree} $\Delta^{tot}(D)$ is the maximum degree of the underlying graph. We write $\Delta^\pm(D)$ for the smallest $m$ such that $D$ has a vertex of outdegree at least~$m$ and a vertex of indegree at least $m$ (these vertices may coincide).
Our main result is the following. 

\begin{restatable}{thm}{oriented}
	\label{thm:oriented}
	Let $T$ be an  oriented tree with $k$ arcs, and  
	let $D$ be a digraph 
	 with  $\delta^0(D)\ge  k/2$ having no oriented $4$-cycles.  
	If $\Delta^\pm(D)> \Delta^{tot}(T)$
	then $D$ contains $T$ as a subgraph.   
\end{restatable}
Observe that  
  the host graph from the example mentioned in the first paragraph can be viewed as a digraph in the standard way (each edge becoming a directed $2$-cycle), and thus  it is necessary to forbid $4$-cycles  
 in Theorem~\ref{thm:oriented}. See Section~\ref{forbid} for a discussion on forbidding a smaller set of oriented $4$-cycles, or a a different family of digraphs.
 
Further,  note that $ \Delta^{tot}(T)$ equals $k$ only if $T$ is a star, and in that case $T$ embeds in $D$ as long as $\Delta^\pm(D)\ge k$ and $\delta^0(D)\ge  k/2$. 
 Hence Theorem~\ref{thm:oriented} has the following corollary, which is the result mentioned in the abstract.
 \begin{restatable}{corollary}{oriented}
	\label{coro:oriented}
	Let $D$ be a digraph 
	 with  $\delta^0(D)\ge  k/2$ and $\Delta^\pm(D)\ge k$ having no oriented $4$-cycles. Then $D$ contains each oriented tree on $k$ arcs as a subgraph.   
\end{restatable}

We provide two more results  for trees that are antidirected or  arborescences.
If the tree we are looking for is antidirected, 
we can replace the minimum semidegree with the \textit{minimum pseudo-semidegree} $\bar\delta^0(D)$ of $D$: this is $0$ if $D$ has no arcs,  and otherwise is defined as the largest integer $d$ such that no vertex $v$ in $D$ has positive indegree $<d$ or positive outdegree $<d$.
Moreover, we  allow  the host $D$ to have directed $4$-cycles. 
\begin{restatable}{thm}{antidirected} 
	\label{thm:antidirected}
	Let $T$ be an antidirected tree with $k$ arcs, and  
	let $D$ be a digraph 
	 with  $\bar\delta^0(D)\ge  k/2$ such that all $4$-cycles of $D$ are directed. 
	If $\Delta^\pm(D)> \Delta^{tot}(T)$
	then $D$ contains $T$ as a subgraph.   
\end{restatable}

In Section~\ref{sec:antitree}, we will deduce from Theorem~\ref{thm:antidirected} the following result:
Each $n$-vertex 
digraph~$D$ with more than $(k-1)n$ arcs whose $4$-cycles  are all directed  contains each  antidirected  tree with $k$ arcs and $\Delta^{tot}(T)\le k/2$.
This is a special case of a conjecture on antidirected trees in digraphs of high arc density from~\cite{ahlrt}. See Section~\ref{sec:antitree} for all details.

If $T$ is an out-arborescence we only need to bound the outdegree of $D$, and for oriented $D$, the bound on the minimum semidegree can be lower than in Theorem~\ref{thm:oriented}. 
\begin{restatable}{thm}{arborescences1}
	\label{thm:arborescences1} 
Let $T$ be a $k$-arc out-arborescence rooted at a vertex of maximum total degree.  
Let $D$ be a digraph having no oriented $4$-cycles with  $\Delta^+(D)> \max\{\Delta^+(T), k/2\}$ and $\delta^+(D)\ge  k/2-\lambda$, where $\lambda=1$ if $D$ is an oriented graph and $\lambda=0$ otherwise. 
 Then $D$ contains $T$ as a subgraph.   
\end{restatable}

A slightly more general statement is given as Proposition~\ref{thm:arborescences} in Section~\ref{sec:arbo}.

The remainder of this paper is organised as follows. After going through some notation in   Section~\ref{sec:notation}, we give the proof of Theorem~\ref{thm:oriented} in Section~\ref{sec:oriented}. 
In Sections~\ref{sec:antidirected}~and~\ref{sec:arborescences} we show how the proof of Theorem~\ref{thm:oriented} can be modified to prove Theorems~\ref{thm:antidirected}~and~\ref{thm:arborescences1}. 
Some open questions and the implication of Theorem~\ref{thm:antidirected} for the conjecture from~\cite{ahlrt} are discussed in Section~\ref{conclusion}.
 
\section{Notation} 
\label{sec:notation}

For a vertex $a$ of a digraph $D$, the \textit{out-neighbourhood} (\textit{in-neighbourhood}; \textit{total neighbourhood}) of $a$, denoted by $N^+(a)$ ($N^-(a)$; $N(a)$) is the set 
of out-neighbours (in-neighbours; out- and in-neighbours) of $a$. %The \textit{neighbourhood} of $a$ is the set $N(a):=N^+(a)\cup N^-(a)$. 
The \textit{outdegree}, \textit{indegree} and \textit{total degree} %(or simply \textit{degree}) of a vertex $a\in D$
of $a$ are  $\deg^+(a):=|N^+(a)|$, $\deg^-(a):=|N^-(a)|$ and $\deg(a):=|N(a)|$. The \textit{minimum outdegree} (\textit{minimum indegree}; \textit{minimum total degree}) of $D$ is $\delta^+(D):=\min_{a\in D}\{\deg^+(a)\}$  ($\delta^-(D):=\min_{a\in D}\{\deg^-(a)\}$; $\delta(D):=\min_{a\in D}\{\deg(a)\}$).  
The \textit{maximum outdegree}  (\textit{maximum indegree}; \textit{maximum total degree}) of $D$ are $\Delta^+(D):=\max_{a\in D}\{\deg^+(a)\}$ ($\Delta^-(D):=\max_{a\in D}\{\deg^-(a)\}$; $\Delta(D):=\max_{a\in D}\{\deg(a)\}$).

Call $D$  {\it antidirected} if each of its vertices has indegree $0$ or outdegree $0$. 
A non-isolated vertex~$a$ of a digraph $D$ is an \textit{out-vertex} (\textit{in-vertex}) if $\deg^-(a)=0$ 
($\deg^+(a)=0$). 
The \textit{distance}  $\dist(a,b)$ of vertices $a,b\in V(D)$ is the length of a
shortest  path of any orientation that joins them. The \textit{diameter} of $D$ is the maximum distance between any pair of its vertices.  If $a$ and $b$ are vertices of  an oriented tree, we denote by $P_{ab}$ the unique path joining them. An \textit{out-arborescence} is an oriented rooted tree in which all edges point away from the root. 
Given an oriented tree $T$, a vertex $u\in V(T)$ is called a \textit{penultimate vertex} if it is not a leaf of $T$ and all its neighbours,
except at most one, are leaves. 

Given an oriented tree $T$ and a digraph $D$, an {\it embedding} of $T$ into $D$ is an injective 
function $f: V(T) \to V(D)$ preserving adjacencies. 
If such an embedding exists, we say that $T$ \textit{embeds} in $D$, or $T$ is a subgraph of $D$.   
For a subset $X\subseteq V(T)$ and an embedding $f$ from $T$ to $D$, we write $f(X)$ to denote
the image of $X$ in $D$.

For convenience, let us define two classes of digraphs, \C and \Ccc \, digraphs. A digraph $D$ is \textit{$\mathcal{C}_4$-free}  if it does not contain any orientation of a $4$-cycle as a subgraph, and it is \textit{$\mathcal{C}^*_4$-free} if any $4$-cycle in $D$ is a directed cycle.

%%%%%%%
%%%%%%%
%%%%%%%
%%%%%%%
%%%%%%%
%%%%%%%
%%%%%%%
%%%%%%%
%%%%%%%
%%%%%%%
%%%%%%%
%%%%%%%
%%%%%%%
%%%%%%%
%%%%%%%
%%%%%%%
%%%%%%%
%%%%%%%
%%%%%%%
%%%%%%%
%%%%%%%
%%%%%%%
%%%%%%%
%%%%%%%
%%%%%%%
%%%%%%%

\section{Proof of Theorem~\ref{thm:oriented}
}
\label{sec:oriented}

In order to be able to easily refer to this proof  
 later, when we prove Theorem~\ref{thm:antidirected}, we will frequently use the weaker property that $D$ is \Ccc \ (instead of using the property that it is \Cc).  
 
 Let $t$ be a  vertex of maximum total degree of $T$.  Without loss of generality we  assume that $\deg^+(t)\ge \deg^-(t)$
	(the other case is analogous).  Let $T_1$ be a maximal subtree of $T$ having diameter at most four 
	that contains $t$ and $N(t)$.	We will show that
	\begin{equation}
		\label{claim:T_1}
		\text{$T_1$ embeds in $D$.} 
	\end{equation}
	For this, let $t'$ be a neighbour of $t$, %with maximum total degree, so 
	and note that $\deg^{tot}(t')\le \lfloor (k+1)/2\rfloor$, with equality if $T$ is a double-star. 
	We  embed   $t$ in a vertex $a\in V(D)$ with $\deg^+(a)\ge  \Delta^+(D)>\deg^{tot}(t)$.  We then embed $t'$ and its neighbours, using the fact that $\delta^0(D)\ge k/2$. Afterwards, we  embed the remaining neighbours of $t$ into unused neighbours of $a$,
	which is possible by our condition on $\Delta^+(D)$	
	and since $D$ is \Cc. If $T=T_1$ is a double-star, then we are done. Otherwise, we proceed by successively embedding $N(u)$ for each non-leaf $u\in V(T_1)$ with an already embedded non-leaf neighbour~$w$. 
	Observe that any already embedded vertex in $T_1$ has  distance  at most two to $w$, 
	because $T_1$ has diameter at most four. 
	So, if $z\in V(T_1)$ is such that $f(z)\in N(f(u))$, then either $z=w$ or $z$ is a neighbour of $w$. Further, there is at most one such neighbour $z$ of $w$,
	because $D$ is \Cc.
Thus, $f(u)$ has at least $k/2 -2$ unused out-neighbours and  at least $k/2-2$ unused in-neighbours, which means we are able to embed $N(u)$. 
This proves~\eqref{claim:T_1}. 
	
	Let$\{T_i\}_{i=1}^r$  be such that $T_r=T$ and each $T_i$ is obtained from $T_{i+1}$
by deleting the leaf neighbours of a penultimate vertex $u_i\in V(T_{i+1})$ with minimum degree. 
	Let $i\in [r]$ be the maximum index such that $T_i$ embeds in $D$. 	For the sake of a contradiction suppose $i<r$. Since $i\ge 1$ by~\eqref{claim:T_1}, the definition of $T_1$ implies that $T_{i+1}$ has diameter at least five. 
	Let $f:V(T_i)\to V(D)$ be an embedding. 
	To simplify our notation, for a vertex  $x \in V(T_i)$ or a subset $X\subseteq V(T_i)$, we denote their image under 
	$f$ by $\td{x}$ and $\td{X}$, respectively, i.e., $\td{x} = f(x)$ and $\td{X} = f(X)$.
	
	Let $u\in V(T_i)$ be the vertex with $\deg_{T_i}(u)<\deg_{T_{i+1}}(u)$, 
	and let $v$ be denote its unique  neighbour in $T_i$. 
	Let $\s\in \{+,-\}$ such that $u\in N^{\s}(v)$. 
	Let $d^+_u$ and $d^-_u$ be the number of leaves of $T_{i+1}$ in $N^+(u)$ and $N^-(u)$, respectively, and let $d_u:=d^+_u+d^-_u=\deg(u)-1$. By symmetry, we can assume that
	\begin{equation}
			\label{du+}
			d^+_u\ge d^-_u. 
		\end{equation}
	
	Let $w\in V(T_i)$ be a penultimate vertex of $T_{i+1}$ with $\dist(v,w)$ maximised. Let
	$W$ be the set of leaves of $T_{i+1}$ in $N(w)$ and let $d_w:=|W|$. By our choice of  $\{T_i\}_{i=1}^r$, we have $d_w\ge d_u\ge 1$, and further,  
	the path $P_{vw}$ from $v$ to $w$ has length at least  two because $T_i$ has diameter at least four. Observe that $u\notin V(P_{vw})$ 
	because $u$ is a penultimate vertex of $T_{i+1}$. 
	Let $v_1$, $v_2$, and $w^*$ denote the second, third, and second-to-last vertices of $P_{vw}$. It may happen that $\{v_1, v_2\} \cap \{w^*, w\} \ne \emptyset$. Define $T' := T_i \setminus \{u\}$.
	
	Set $Q:=N^\diamond(\td{v})\setminus f(V(T'))$ and $q:=|Q|$. Note that $q\ge 1$ as $\td{u}\in Q$. 	
	Consider any $a\in Q$. If there exist disjoint sets $S_a^+\subseteq N^+(a)\setminus f(V(T'))$ and $S_a^-\subseteq N^-(a)\setminus f(V(T'))$
	with $|S_a^+|\ge d_u^+$ and $|S_a^-|\ge d_u^-$, then it is straightforward to embed $T_{i+1}$, 
	 contradicting the maximality of $i$. 
	Thus, such sets $S_a^+$ and $S_a^-$ do not exist, and therefore at least one of $|N^+(a)\setminus f(V(T'))|$ or $|N^-(a)\setminus f(V(T'))|$ is smaller than $d_u$.
	 If $|N^+(a)\setminus f(V(T'))|<d_u$ then set $N_{ a}:=N^+(a)$. Otherwise, note that the embedding of $T_{i+1}$ failed because $T_{i+1}$ has a leaf  in $N^-(u)$, and set  $N_{a}:=N^-(a)$. Observe that for distinct $a_1,a_2\in Q$, either $N_{a_1}$ and $N_{a_2}$ are both out-neighbourhoods or both in-neighbourhoods. This means that for all $a \in Q$, we consistently consider either the out-neighbourhoods or the in-neighbourhoods. 
	 In either case, we have
	\begin{equation}
		\label{eq:neighbours-Q-oriented}
		|N_{a}\cap f(V(T'))|\ge   k/2 -d_u+1 \textrm{\, for each $a\in Q$.} 
	\end{equation}
	
	Setting $$R:=f(V(T')\setminus (W\cup \{v\}),$$  we have   
	\begin{equation}\label{RRR}
		2\le |R|\le k+1-d_u-1-d_w-1=k-d_u-d_w-1.
	\end{equation}
	Let $R'\subseteq R\subseteq V(D)$ be the set of vertices that do not belong to $N_a$ for any $a\in Q$.  
	Observe that $\td{v_2}\in R'$, because $D$ is \Cc, and so $|R'|\ge 1$. 
	Then, 
		\begin{equation}\label{QQQ}
	q(k/2-d_u+1)\le \sum_{a\in Q}|N_{ a}\cap f(V(T'))| \le  |R|-|R'|+\min\{q,d_w\}+q \le  |R|-|R'|+2q,
	 \end{equation}
	where the second inequality holds because $f(V(T'))=R\cup \hat W\cup\{\hat v\}$ and because  
	\begin{itemize}
		\item \label{item:b3} $\sum\limits_{a\in Q}|N_{ a}\cap R|\le  |R|-|R'|$, as each vertex $b\in R-R'$ belongs to exactly one set $N_{a}$, and
\item \label{item:b2} $\sum\limits_{a\in Q}|N_{a}\cap \td{W}|\le \min\{d_w,q\}$ as  each $a\in Q$ 
		has at most one neighbour $b\in \td{W}$ and vice versa. 
	\end{itemize}
	We also used the the fact that  $D$ is \Ccc. Indeed, if there exist distinct $a_1,a_2\in Q$ and $b\in R-R'$ such that 
	$b\in N_{a_1}\cap N_{a_2}$,  then the subset $\{a_1,a_2,\td v, b\}$ would induce an oriented $4$-cycle (in $D$) that is not directed. 
	Moreover, if the second inequality does not hold, then either there exist $a_1,a_2\in Q$ and $w\in W$ such that the set $\{\td v,a_1,a_2,\td w\}$ induces a $4$-cycle in $D$ that is not directed,  or there exist $a\in Q$ and $w_1,w_2\in W$ such that the set $\{a,\td w_1,\td w_2,\td w\}$ induces a $4$-cycle which is not directed, leading to a contradiction in both cases. 
	 
	Simplifying~\eqref{QQQ}, and using~\eqref{RRR} as well as the facts that $|R'|\ge 1$ and $d_w\ge d_u$, we obtain
	\begin{equation}\label{aux}
		q(k/2-d_u-1) \le  k-d_u-d_w-1-|R'| \le  k-2d_u-2 =2(k/2-d_u-1).
	\end{equation} 
	Thus, since $d_u+1=\deg(u)< k/2$ (as $\deg(u)\le \deg(t)$ and $T$ has diameter at least five), 
	we have
	\begin{equation}
		\label{qle2}
		q\le 2. 
	\end{equation} 
Furthermore, if $q=2$, all inequalities in~\eqref{QQQ} and~\eqref{aux}
 must hold with equality. Therefore, if $q=2$, then
		\begin{enumerate}[label=$(a\arabic*)$]
			\item \label{item:a1} $q\le d_w=d_u$, 
			\item \label{item:a2} $|N_{\hat u}\cap \td{W}|=1$, 
			\item \label{item:a3} for each neighbour $h$ of $u$ which is a leaf of $T_{i+1}$, there is an embedding $g$ of the tree $T^*=T_{i+1}-h$ with  $g(t)=f(t)$ for each $t\in T_i$.
		\end{enumerate}
		For~\ref{item:a3}, observe that   $|N_{\hat u}\cap f(V(T'))|=k/2-d_u+1$ implies that 
	\[|N^+(\hat u)\setminus f(V(T'))|, |N^-(\hat u)\setminus f(V(T'))| \ge d_u-1. 	\]
	
	Set $N_{\td{v}}:=N^\s(\td{v})$. 
	Since $\delta^0(D)\ge k/2$ and since $D$ is \Ccc, we have 
		\begin{equation}\label{NuNvNuv}
			|N_{\td{u}}\cap \im f|\ge k/2-d_u+1, \quad%\textrm{ and } \quad 
			|N_{\td{v}}\cap \im f|\ge k/2-q+1\quad\textrm{ and } \quad |N_{\td{u}}\cap N_{\td{v}}|\le 1.
		\end{equation}
	 
	The remainder of the proof  splits into two cases. 		\smallskip	
		\\
		\textbf{Case A: $\boldsymbol{q=2}$.}  By~\ref{item:a2}, 
		there exists $w_1\in W$ with $\td{w_1}\in N_{\td{u}}$. 
		Assume $w_1\in N^+(w)$
		(otherwise change $N^+(\td{w})$ to $N^-(\td{w})$ in what follows). 
		Set $N_{\td{w}}:=N^+(\td{w})$. 	
		We have,
		\begin{equation}\label{vwnot}
		\text{$\hat w\notin N_{\td v}$  and $\td v\notin N_{\td w}$,} 
		\end{equation}
		 since $D$ is \Ccc.  
			We claim that
		\begin{equation}
			\label{eqq:number-neighbours}
			 |N_{\td{w}}\cap \im f|\ge k/2.
		\end{equation}
		Indeed, if this inequality does not hold, consider the embedding of $T^*$ from~\ref{item:a3}, where $h$ is chosen as follows. If $N_{\hat u}=N^+(u)$ we let $h$ be an out-neighbour of $u$ that is a leaf in $T_{i+1}$ (such a vertex $h$ exists because of~\eqref{du+}). If $N_{\hat u}=N^-(u)$, then our choice of $N_{\hat u}$ guarantees the existence of an in-neighbour $h$ of $u$ that is a leaf in $T_{i+1}$.
		Now, modify the embedding of $T^*$ as follows: embed $w_1$ in an unused vertex of $N_{\td{w}}$, and embed
		$h$ in $\hat w_1$. This is an embedding of all of $T_{i+1}$, a contradiction, so~\eqref{eqq:number-neighbours} is proved.

		Since $D$ is \Ccc, we have $|N_{\td{w}}\cap N_{\td{z}}|\le 1$, for any $z\in \{u,v\}$. 
		So,  \eqref{NuNvNuv},~\eqref{vwnot}  and~\eqref{eqq:number-neighbours} give 
		\begin{equation}
			\label{eq:caseA}
			k-d_u+1\ge |\im f|\ge  |\{v,w\}|+|N_{\td{v}}\cap \im f| + |N_{\td{w}}\cap \im f| -1\ge k
		\end{equation}
		where all inequalities must hold as equalities. Thus $d_u=1$, a contradiction to~\ref{item:a1}. 
\smallskip \\ 
	\textbf{Case B: \bf $\boldsymbol{q=1}$.} 
			Let $Y\subseteq V(T_i)$ such that $\td{Y}=f(V(T_i))\setminus (N_{\td{u}}\cup N_{\td{v}})$. Consider the subpath $w^*,w,w_1$ of $T_i$. Suppose, for contradiction, that $Y=\emptyset$, so $\td w^*,\td w,\td w_1\in N_{\td{u}}\cup N_{\td{v}}$. Without 
				loss of generality, suppose $\td w^*\in N_{\td{u}}$. Since $D$ is \Ccc, we must have $\td w_1\in N_{\td{v}}\setminus N_{\td{u}}$. 
				Then, either $\td w\in N_{\td{u}}$ or $\td w\in N_{\td{v}}$. In the former case, the set $\{\td v,\td u,\td w,\td w_1\}$ induces
				a $4$-cycle that is not directed, whereas in the later case, the set $\{\td v,\td u,\td w,\td w^*\}$ induces a $4$-cycle other than
				the directed one; both cases contradict that $D$ is \Ccc. Therefore, we must have $Y\ne \emptyset$.   
		
		We have
	\begin{equation}\label{YYY}
		k-d_u+1\ge |V(T_i)|=|\im f|\ge |N_{\td{u}}\cap \im f|+| N_{\td{v}}\cap \im f|-1+|\td{Y}|\ge k-d_u+|\td{Y}|.\end{equation}
	Thus
	\begin{equation}\label{1Yq}
		|Y|= 1.\end{equation}

		Let $Y=:\{y\}$. Then all  inequalities in~\eqref{NuNvNuv} and~\eqref{YYY} hold as equalities. In particular, there is exactly one vertex $x\in V(T_i)$ with $\hat x \in N_{\td{u}}\cap N_{\td{v}}\cap \im f$ and  $|V(T_i)|=k-d_u+1$, i.e., \begin{equation}\label{i=r-1}i=r-1.\end{equation}

		Let $J$ be the subgraph of $T_i$ induced by~{$V(T_i)\setminus \{u,v,y\}$}. We claim that
		\begin{equation}
			\textrm{$J$ contains no oriented path of length two, and $x$ is an isolated vertex in $J$.}
			\label{eq:J}
		\end{equation}
		Indeed,
		any oriented path $z_1z_2z_3$ of length two in $J$
		either contains an arc  $z_iz_{i+1}$  such that 
		$\td{z_i}\in N_{\td{u}}$ and $\td{z}_{i+1}\in N_{\td{v}}$, or has the property that $\td{z_1},\td{z_3}$ are either both in $N_{\td{u}}$ or both in  $N_{\td{v}}$. Either option contradicts the fact that $D$ is \Ccc. 
		Moreover, if $x$ has a neighbour in $J$, we found a $4$-cycle that is not directed, again
		a contradiction. This proves~\eqref{eq:J}.
		
		Next, we claim that
		\begin{equation}
		d_w=1.
			\label{Wis1}
		\end{equation}
		Suppose otherwise. Then there exist distinct vertices  $w_1,w_2\in W$. 
		Note that if $y\ne w$, then $w$ and any two vertices in $\{w_1,w_2,w^*\}\setminus \{y\}$ together
		induce an oriented path of length two in $J$,
		contradicting~\eqref{eq:J}. On the other hand, 
		if $y=w$, then  two of $\{\td{w_1},\td{w_2},\td{w}^*\}$ must belong either to  $N_{\td{u}}$ or to $N_{\td{v}}$,  a contradiction to $D$ being \Ccc. This proves~\eqref{Wis1}.
		
	By~\eqref{Wis1}, the set $W$ contains only one element, say $w_1$.
		Then, by~\eqref{eq:J},
		\begin{equation}
			\label{eq:y}
			y\in \{w^*,w,w_1\}. 
		\end{equation}
		Moreover, by~\eqref{Wis1} and our choice of $w$, 	we know that  $d_u=1$. %By symmetry we can assume that $d_u^+=1$  and $d_u^-=0$. 
		So, by~\eqref{du+},  we have $d_u^+=1$  and $d_u^-=0$. By our choice of $N_{\hat u}$ this means that $N_{\hat u}=N^+(\td u)$.

		We claim that
		\begin{equation}
			\label{eq:xinNv}
			x\in N^\s(v). 
		\end{equation}
		Suppose otherwise. Then the path  $P_{vx}$ from $x$ to $v$ has length at least $2$. %\A{acá cambié ligeramente para tener siempre \Ccc, el caso de que $P_{vx}$ tenía longitud 3 se puede incorporar en el caso en que tenía longitud al menos 4} 
		Given that $D$ is \Ccc\, and $\hat x\in N_{\hat u}\cap N_{\hat v}$, and by~\eqref{eq:J}, we know  $P_{vx}$ cannot have length $2$, so it must have length at least $3$. By~\eqref{eq:J}, this implies that in $T_i$, the only neighbour of $x$ is $y$. 
		Thus, by~\eqref{eq:y}, we have $x=w_1$, $y=w$, and then either $w^*\in N_{\hat u}$ or $w^*\in N_{\hat v}$, implying the existence of an oriented $4$-cycle in $D$ that is not directed. This contradicts the fact that $D$ is \Ccc, 
		and thus proves~\eqref{eq:xinNv}.  
		
		By~\eqref{eq:J},  $x$ is isolated  in $J$,  so  if $x=v_1$ then $y=v_2\ne w_1$. Letting $z$ be a
			neighbour of $y$ in $T_{i}$, we see that either $\{\hat u,\hat x,\hat y,\hat z\}$ or $\{\hat v,\hat x,\hat y,\hat z\}$ span a $4$-cycle, contradicting that $D$ is \Ccc.  Thus,
			\begin{equation}
			\label{eq:xneqv1}
			x\ne v_1. 
		\end{equation}
		So $x\notin V(P_{uw})$. Observe that if $|N^+(\td{x})\setminus \im f|\ge 1=d_u$, then 
			we can embed $T_{i+1}$ by swapping the images of $u$ and $x$, and embedding $u$'s unembedded neighbour into an
			unused vertex in $N^+(\td{x})\setminus \im f$, a contradiction. 
			Therefore, $N^+(\td{x})\subseteq  \im f$, and hence 
			$|N^+(\td{x})\cap \im f|\ge  k/2$. 
			Because $d_u=1$ and by~\eqref{i=r-1}, we have $|V(T_i)|=k$, and so, using the fact that $D$ is \Ccc, and the inequalities from~\eqref{NuNvNuv}, we see that
			\[k= |\im f|\ge |N_{\td{u}}\cap \im f|+|N_{\td{v}}\cap \im f|+|N^+(\td{x})\cap \im f|-3\ge 3k/2-3,\]
			implying that $|V(T_i)|= k\le 6$. 
			Thus, $V(T_i)=\{u,v,v_1,w,w_1,x\}$. 
			As $D$ is \Cc,
			we have $\td v_1,\td w\notin N_{\td u}\cup N^+({\td x})$, 
				and further, $\td w_1\notin N_{\td u}\cap N^+({\td x})$, implying that either $|N_{\td u}\cap \im f|< 3$ or $|N^+({\td x})\cap \im f|<3=k/2$, a contradiction.

\section{Proof of Theorem~\ref{thm:antidirected}: Embedding Antidirected Trees} 
\label{sec:antidirected}
The proof   of Theorem~\ref{thm:antidirected} follows along the lines of the proof of Theorem~\ref{thm:oriented}, and we will limit ourselves to pointing out the differences and how we resolve them. Note that we now only have to embed antidirected trees, but have to cope with the weaker condition that $D$ is \Ccc  \ instead of \Cc , which means that $D$ is allowed to contain directed $4$-cycles. Further, we only have the lower bound of $k/2$ on the minimum pseudo-semidegree $\bar\delta^0(D)$, instead of the minimum semidegree $\delta^0(D)$. 

Let us first argue why allowing directed $4$-cycles is not a problem when embedding antidirected trees. For this, we check all places in the proof of Theorem~\ref{thm:oriented} where we actually use the property that the host digraph $D$ is \Cc, and show that, when $T$ is antidirected,  replacing this condition with $D$ being \Ccc\, allows the same argument to hold.  
\begin{itemize}
	\item When embedding $N^+(t)\cup N^-(t')$, we observe that $|N^+(\td t)\cap N^-(\td t')|\le 1$ because $D$ is \Ccc. By our condition on $\Delta^+(D)$, it is feasible to first embed $N^-(t')\cup \{t'\}$ so that $t$ is mapped to a vertex $a$ with $\deg^+(a)>\deg^+(t)$, followed by the remaining unembedded neighbours of $t$. 
	\item When embedding $N^+(u)$ for a non-leaf vertex $u$ of $T_1$ with an already embedded neighbour $w$, and assuming without loss of generality that $u$ is an out-vertex, if $z$ is a vertex of $T_1$ with $\td z\in N^+(\td u)$, then, since $D$ is \Ccc, it follows that either $z=w$ or $z\in N^-(w)$. Thus, at most two vertices of $N^+(\td u)$ are already used, ensuring that $N^+(\td u)$ has enough unused vertices to embed the unembedded out-neighbours of $u$.
	\item  When considering the subset $R'$, we again have $v_2\in R'\ne \emptyset$ as desired, because $D$ is \Ccc.  
	\item In the last paragraph of the proof, when $k=6$ and 
	$V(T_i)=\{u,v,v_1,w,w_1,x\}$, given that $D$ is \Ccc, we have 
	$\td v_1,\td w\notin N_{\td u}\cup N^+({\td x})$
	and $\td w_1\notin N_{\td u}\cap N^+({\td x})$. This implies that either $|N_{\td u}\cap \im f|< 3$ or $|N^+({\td x})\cap \im f|<3=k/2$, a contradiction. 
	\end{itemize}
Next, we explain why  $\delta^0(D)$   can be replaced with $\bar\delta^0(D)$. Specifically, we will verify the following. 
\begin{equation}
	\label{eq:anti-second_claim}
	\begin{minipage}{0.8\textwidth}
		\textrm{
			In the proof of Theorem~\ref{thm:oriented}, whenever the tree $T$ is antidirected, it  suffices to consider, for every vertex $v$ in $D$, either its out-neighbourhood (if $d^+(v)>0$) or its in-neighbourhood (if $d^-(v)>0$). 
			}
	\end{minipage}
\end{equation}

When embedding the subtree $T_1$, we begin by embedding the double star with centers $t$ and $t'$. Without loss of generality, assume that $t$ is an out-vertex; otherwise interchange $+$ and $-$ in the following argument. Once the arc $tt'$ has been embedded into the arc $\td t \td t'$, where $\deg^+(\td t)>\deg^+(t)$, we consider only the sets $N^+(\td t)$ and $N^-(\td t)$, ensuring that~\eqref{eq:anti-second_claim} holds up to this step. 
Next, when embedding the neighbours of a non-leaf vertex $u$ of $T_1$ that has an already embedded neighbour $w$, it suffices to consider either $N^+(\td u)$ or $N^-(\td u)$, depending on whether $u$ is an out-vertex or an in-vertex. Moreover, the required neighbourhood exists: since $u$ and $w$ are already embedded, $\td u$ has outdegree or indegree at least $k/2$ because $\bar\delta^0(D)\ge k/2$.  Thus,~\eqref{eq:anti-second_claim} holds for this step as well. 
 
 Let us now examine the structures that are considered in the subsequent argument.  Without loss of generality, assume that $u$ is an out-vertex, and therefore $v$ is an in-vertex. First, note that the set $Q=N^-(\td v)\setminus f(V(T'))$ is well-defined and non-empty, as $v$ is an already embedded in-vertex and $\td u\in Q$. Moreover, we have $N_{\td v}=N^-(\td v)$.  Second,  for each $a\in Q$, we have the set $N_a=N^+(a)$, that is non-empty because $\td v\in N_a$ as $u$ is an out-vertex, and further, $\deg^+(a)\ge k/2$ because $a\in Q\subseteq N^-(\td v)$. Next, in Case A, we have $N_{\td w}=N^+(\td w)$, where $\deg^+(\td w)\ge k/2$ as $w$ and its out-neighbours  are already embedded. 
Finally, in Case B, we additionally consider the set $N^+(\td x)$. Since $x$ is an already embedded in-neighbour of $v$, we have $\deg^+(\td x)\ge k/2$. 
In conclusion, each of the sets $N_{\td u}$, $N_{\td v}$, $N_{\td w}$ and $N^+(\td x)$, which are considered in addition to the subgraph $\im f$ of~$D$, is well-defined and has cardinality at least $k/2$ by the minimum pseudo-semidegree condition of~$D$, ensuring that~\eqref{eq:anti-second_claim} holds. This concludes the proof.

\section{Proof of Theorem~\ref{thm:arborescences1}: Embedding Arborescences}  
\label{sec:arborescences}
\label{sec:arbo}
We will prove the following result which implies Theorem~\ref{thm:arborescences1}.
\begin{restatable}{prop}{arborescences}
	\label{thm:arborescences} 
Let $T$ be a $k$-arc out-arborescence rooted at a vertex of maximum total degree.  
Let $D$ be a digraph having no oriented $4$-cycles with  $\Delta^+(D)> \Delta^+(T)$. If 
\begin{itemize}
\item $\delta^+(D)\ge  k/2-1$, $\Delta^+(D)\ge k/2$ and $D$ is an oriented graph, or
 \item $\delta^+(D)\ge  k/2$,
\end{itemize}
 then $D$ contains $T$ as a subgraph.   
\end{restatable}

	The proof of Proposition~\ref{thm:arborescences} proceeds similarly to the proof of Theorem~\ref{thm:oriented}. Again, we just highlight  the key differences.   The following notation will be useful: For an out-arborescence $T$ and $u\in V(T)$ which is not its root, the \textit{parent}~$p_u$ of $u$ is the unique in-neighbour of $u$. 
	
	Let us first establish that $T_1$, rooted at $t$, embeds in $D$. We consider both cases at once: the case 
	when $D$ is oriented and the case when $D$ is an arbitrary digraph. 
	We start by embedding the double star with centers $t$ and $t'$. After embedding the arc $tt'$ into the arc $\td t \td t'$, we consider only the sets $N^+(\td t)$ and $N^+(\td t')$ to embed the remaining out-neighbours of $t$ and $t'$, respectively. This is feasible because $\deg^+(\td t)>\deg^+(t)$ and $\deg^+(\td t')\ge \deg^+(t')$.  Additionally, if $D$ is an oriented graph, we require that $\deg^+(\td t)\ge k/2$, which is ensured by our hypothesis. Next, when embedding the out-neighbours of a non-leaf vertex $u$ of $T_1$ whose parent has already been embedded, we only  consider the set $N^+(\td u)$. This set contains at least $k/2-2$ unused vertices:  if $D$ is oriented, at most one vertex in $N^+(\td u)$ has already been used, and otherwise, at most two vertices have been used. Thus, it is feasible to embed the out-neighbours of $u$ into unused vertices in $N^+(\td u)$, ensuring that $T_1$ embeds in $D$. 
	
	Let $i$ be the maximum index such that $T_i$ embeds in $D$, with the additional requirement that  if~$D$ is oriented, then $\deg^+(\td t)\ge k/2$. From the previous argument, we know that $i\ge 1$. Now, suppose for the sake of a contradiction that $T_i\ne T$. 
	
	 First, suppose $D$ is not an oriented graph.
	We proceed along the lines of the proof of Theorem~\ref{thm:oriented}.
	Define $Q=N^+(\td v)\setminus f(V(T'))$, which is non-empty as both $u$ and $v$ are already embedded, with $\td u\in Q$. 
	We set $N_{\td v}=N^+(\td v)$ and $N_a=N^+(a)$ for each $a\in Q$. Each of these sets is well-defined and contains at least $k/2$ vertices due to the minimum outdegree condition of $D$. Additionally, we  consider the sets $N_{\td w}=N^+(\td w)$, and $N^+(\td x)$ for Case B. These sets are also well-defined, and each contains at least $k/2$ vertices since $\delta^+(D)\ge k/2$. It is straightforward to verify that the sets $N_{\td u}$, $N_{\td v}$, $N_{\td w}$, $N^+(\td x)$ satisfy the conditions required in the proof of Theorem~\ref{thm:oriented}. This leads to a contradiction, which completes the proof in the case when $D$ is not an oriented graph.
	
	Now, suppose $D$ is an oriented graph. 
	Roughly speaking, the reason why we can decrease the out-degree requirement by one in this case is that, 
	for any vertex $z$ that has already been embedded, 
	$\td p_z$ does not belong to the out-neighbourhood $N^+(\td z)$. This contrasts with the non-oriented case, where it is possible that  
	$\td p_z\in N^+(\td z)$, regardless that $p_z$ is, indeed, an in-neighbour of $z$.  
	
	Note that, among vertices in $\{u,v,w,x\}$, only $w$ may be the root $t$ of $T$. 
	Define $Q, N_{\td v}, N_a, N_{\td w}$ as in the case when $D$ is non-oriented. Since $N_a\subseteq N^+(v)$ and $D$ is oriented, 
	we have $\td v\notin N_a$ for any $a\in Q$. Given this distinction, and following the approach in the proof of 
	Theorem~\ref{thm:oriented}, it is straightforward to verify that $q\le 2$, and furthermore, if $q=2$ then~\ref{item:a1}-\ref{item:a3} remain valid. 
	
	Let us argue that $q\neq 2$. Indeed, otherwise
	by~\ref{item:a3} we have $N_{\td w}\subseteq \im f$, and so, 
	\[|\im f|\ge
	\begin{cases}
		 |\{v,w,p_v,p_w\}|+|N_{\td v}\cap \im f|+|N_{\td w}\cap \im f|-1\ge k, &\textrm{ if $w\ne t$,} \\
		|\{v,w,p_v\}|+|N_{\td v}\cap \im f|+|N_{\td w}\cap \im f|-1\ge k, &\textrm{ if $w=t$. } \\
	\end{cases}
	\]
	Since $|\im f|\le k-d_u+1$, it follows that $d_u=1$, which contradicts~\ref{item:a1}. 
	
	So $q=1$. Define $Y:=f(V(T_i))\setminus (N_{\td u}\cup N_{\td v}\cup \{\td p_u, \td p_v\})$, where $p_u=v$. 
	Similarly as in the proof of Theorem~\ref{thm:oriented}, we see that $|Y|=1$ 
	and   that there is a unique vertex
	$\td x\in (N_{\td u}\cup \{\td p_u\})\cap (N_{\td v}\cup \{\td p_v\})$. Define~$J$ as in the proof of Theorem~\ref{thm:oriented}. 
	It is straightforward to verify   statements~\eqref{eq:J}--\eqref{eq:xneqv1}. Indeed, these arguments 
	rely solely on the structure of $T$ and the property that $D$ is \Ccc, and thus also \Cc; the minimum out-degree condition of $D$ is not used.  
	Additionally, 
	we have $N_{\td x}=N^+(\td x)\subseteq \im f$. 
	Note that, since $\td x$ and $\td u$ share $\td v$ as a common in-neighbour, it follows that
	$N_{\td x}\cap N_{\td u}=\emptyset$. Furthermore, if $p_v\in N_{\td u}\cup N_{\td x}$, then  $N_{\td x}\cap N_{\td v}=\emptyset$, 
	or $N_{\td u}\cap N_{\td v}=\emptyset$.  Thus,
	\[|\im f|\ge
		\begin{cases}
	|\{\td v,\td p_v\}|+|N_{\td{u}}\cap \im f|+|N_{\td{v}}\cap \im f|+|N_{\td{x}}\cap \im f|-2\ge 3k/2-3, & \textrm{if $p_v\notin N_{\td u}\cup N_{\td x}$, }\\
	|\{\td v\}|+|N_{\td{u}}\cap \im f|+|N_{\td{v}}\cap \im f|+|N_{\td{x}}\cap \im f|-1\ge 3k/2-3, & \textrm{if $p_v\in N_{\td u}\cup N_{\td x}$. }\\
		\end{cases}
	\]
	Since $|\im f|\le k$, we conclude that 
	$k=6$ and $V(T_i)=\{u,v,v_1,w,w_1,x\}$, which contradicts the fact that $T$ is an arborescence rooted at a vertex of maximum total degree. 
	This  completes the proof.

\section{Final remarks}\label{conclusion}

\subsection{Implications for  a conjecture on antidirected trees}
\label{sec:antitree}
In 2013,
Addario-Berry, Havet, Linhares Sales, Reed and Thomassé~\cite{ahlrt} proposed the following conjecture, which by previous work of Burr~\cite{burr} would be optimal.
\begin{conjecture}[Addario-Berry et al.~\cite{ahlrt}]\label{antitreeconj}
Every  
digraph $D$ with more than $(k-1)|V(D)|$ arcs contains each  antidirected  tree with $k$ arcs.
\end{conjecture}
 This conjecture was confirmed in~\cite{ahlrt}
for stars and double-stars, and an asymptotic version for large bounded degree trees was proved in~\cite{camila}. 
Very recently, the present authors established in~\cite{stein-trujillo} two special cases of  Conjecture~\ref{antitreeconj}. First, it holds whenever the antidirected tree is 
a caterpillar. Second, it holds whenever the host digraph $D$  does not contain  any of three specific orientations of
$K_{2,\lceil k/12\rceil}$, namely the ones where each of the two vertices in the partition class of size two has either outdegree~$0$ or indegree $0$. 

Theorem~\ref{thm:antidirected} in the present paper has the following corollary (for a sketch of a proof see below), which states that 
Conjecture~\ref{antitreeconj} holds if all $4$-cycles of $D$ are directed (this includes the case that $D$ has no oriented $4$-cycles at all) and $\Delta^{tot}(T)\le k/2$:
\begin{corollary}\label{coroanti}
Every  
digraph $D$ with more than $(k-1)|V(D)|$ arcs  whose $4$-cycles  are all directed  contains each  antidirected  tree $T$ with $k$ arcs and $\Delta^{tot}(T)\le k/2$.
\end{corollary}
Note that for $k\ge 13$, this corollary also follows from the results of~\cite{stein-trujillo}, as then each of the forbidden  orientations of $K_{2,\lceil k/12\rceil}$ contains a $4$-cycle that is not directed. 

We deduce Corollary~\ref{coroanti} from  Theorem~\ref{thm:antidirected} as follows.
Given an antidirected tree $T$ with $k$ arcs and a digraph $D$ with  more than $(k-1)|V(D)|$ arcs and whose only $4$-cycles are directed,   we  use Lemma~9 from~\cite{antipaths},
 which says that every $n$-vertex digraph with  more than $(k-1)n$ arcs has a subdigraph of minimum pseudo-semidegree at least $k/2$. Applying this lemma to $D$, we obtain a subdigraph $D'$ of~$D$ with $\bar\delta^0(D')\ge k/2$ whose only $4$-cycles are directed. Then  Theorem~\ref{thm:antidirected} provides an embedding of~$T$ in $D'$, and thus in $D$.

\subsection{No forbidden subgraphs}\label{sec:otherdegree}
Let us now discuss alternatives to forbidding all oriented  $4$-cycles. We first describe the situation for graphs. As we saw in the beginning of the introduction, a graph of minimum degree at least $k/2$ and maximum degree at least $k$ does not necessarily contain all trees on $k$ edges. However, this  changes if we raise the bound on the minimum degree to almost $k$: 
 In~\cite{havet},  it was shown that there exists a $\gamma>0$ such that every graph $G$ with  $\delta (G)\ge (1-\gamma)k$ and $\Delta(G)\ge k$ contains all trees on $k$ edges.  
 Havet, Reed, Stein and Wood~\cite{havet} conjectured that the bound on the minimum degree from the previous sentence can be lowered to $\lfloor 2k/3\rfloor$, and this conjecture is supported by evidence found in~\cite{bps,  brucemaya1, brucemaya2}. 
 The bound $\lfloor 2k/3\rfloor$ would be  tight, which can be seen by a small modification of the example from the introduction: Take two copies of the clique  on $\lfloor 2k/3\rfloor -1$ vertices and add a universal vertex, and note that the resulting  graph does not contain any $k$-edge tree having a vertex with three   equal-sized components attached to it.

In~\cite{maya-digraphs} the following digraph analogue was suggested:
\begin{problem}[\cite{maya-digraphs}] 
	\label{prob:f(k)}
	Determine the smallest $f(k)$ such that for every natural $k$, every  digraph of minimum semidegree 
	exceeding $f(k)$ that has a vertex $v$ with $d^+(v),d^-(v)\ge k$ contains each oriented $k$-edge tree. 
\end{problem}
Since we can embed any  oriented $k$-edge tree greedily in any digraph of minimum semidegree at least $k$, we know that $f(k)\le k$. 
In the forthcoming paper~\cite{giovannegeorgemaya}, Kontogeorgiou, Santos and the first author show that asymptotically, $f(k)$ can be taken as $2k/3$ if the tree is large and has bounded degree. 
This value would be best possible  because of the example from the first paragraph of this section, where we view the host  as a digraph and the tree as an oriented tree $T$.  

Note that $T$ does not need to have large degrees: it can be chosen to have total maximum degree~$3$, and this does not affect the example. 
For oriented paths, however, the example fails to work. It is possible that  for oriented paths we can leave the minimum semidegree bound at $k/2$, and omit the bound on the maximum degree. Indeed, a conjecture from~\cite{survey} states that every digraph~$D$ with $\delta^0(D)> k/2$ contains each $k$-edge oriented path. For progress on this conjecture, see~\cite{antipaths,chen2024, skokanTyomkyn}.

\subsection{Forbidding other classes of digraphs}\label{forbid}

We do not know whether in Theorem~\ref{thm:oriented}, it is necessary to forbid all  orientations of the $4$-cycle. In our proof of Theorem~\ref{thm:oriented}, most of the arguments  only require the host $D$ to be \Ccc\ (thus allowing for 
the presence of directed $4$-cycles). Perhaps the few places where we relied on forbidding directed $4$-cycles can be handled differently, and then, directed $4$-cycles can be allowed.
Possibly also other types of oriented $4$-cycles can be allowed. 
We pose the following problem. 
\begin{problem}\label{probo}
	Determine a minimal family $\mathcal C_4'$ of oriented $4$-cycles such that  each $\mathcal C_4'$-free digraph $D$ of minimum 
	semidegree at least $k/2$ and with $\Delta^\pm(D)\ge k$ contain each  $k$-edge oriented tree. 
\end{problem} 

Note that at least one type of oriented $4$-cycles has to be forbidden, as the example from the introduction showed. In the host digraph from this example all four different types of oriented $4$-cycles appear. Interestingly,  in the following example only two types of oriented $4$-cycles appear: Take a directed cycle of any length and substitute each vertex with an independent set of $k/2$ vertices. Arcs are inherited from the original directed cycle. It is not hard to see that this oriented graph $O$ fails to contain any antidirected tree whose partition classes have different sizes. The oriented $4$-cycles contained in~$O$ are exactly  those that  contain no directed path with $3$ arcs. So, the family  $\mathcal C_4'$ from Problem~\ref{probo} needs to contain at least one of these two oriented  cycles.

Another possibility of weakening the conditions on $D$ may be to forbid complete bipartite graphs. In the spirit of~\cite{stein-trujillo}, the family $\mathcal K_{2,s}$   of 
	 all orientations of $K_{2,s}$ seems a natural candidate.

\begin{question}
	For which $s\ge 3$ does each $\mathcal K_{2,s}$-free digraph $D$ of minimum 
	semidegree at least $k/2$  and with $\Delta^\pm(D)\ge k$ contain each  $k$-edge oriented tree?
\end{question} 

As mentioned in the introduction, Brandt and Dobson~\cite{erdos-sos-girth5} showed that every graph of girth at least five with $\delta(G) \ge k/2$ and $\Delta(G) \ge \Delta(T)$ contains every $k$-edge tree. Dobson~\cite{dobson} further conjectured that any graph of girth at least $2\ell + 1$ with $\delta(G) \ge k/\ell$ and $\delta(G)\ge \Delta(T)$ contains each $k$-edge tree, a result confirmed by Jiang~\cite{Jiang01} after prior work by Haxell and {\L}uczak~\cite{hax-luc}. Note that the condition $\delta(G)\ge \Delta(T)$ (instead of $\Delta(G)\ge \Delta(T)$) is indeed necessary, which can be seen by considering a balanced double star and a suitable host graph. Graphs with other kinds of expansion properties have also been studied, see for instance~\cite{sudakov}. 

In the digraph setting, the following question would be a natural analogue of Jiang's result: 
\begin{question}[\cite{maya_LAwomen}]
	Is it true that every oriented graph $D$ 
	of girth at least $2\ell + 1$ with $\delta^0(D) \ge \max\{k/\ell, \Delta^{tot}(T)\}$
	contains every $k$-edge oriented tree?
\end{question}

\small\bibliographystyle{abbrvnat} \bibliography{biblio}

\end{document}